\documentclass[12
pt,twoside,a4paper]{amsart}
\usepackage{amssymb}
\usepackage[latin1]{inputenc}
\usepackage[T1]{fontenc}
\usepackage{times}

\def\hpq0{h^{p,q}_{\leq 0}}
\def\Hpq0{\H_{\leq 0}^{p,q}}

\def\dbar{\bar\partial}
\def\ddbar{\partial\dbar}

\def\C{{\mathbb C}}

\def\H{{\mathcal H}}

\def\Re{{\rm Re\,  }}

\def\be{\begin{equation}}
\def\ee{\end{equation}}

\newtheorem{thm}{Theorem}[section]

\theoremstyle{definition}

\theoremstyle{remark}

\newtheorem{preremark}{Remark}
\newtheorem{preex}{Example}

\numberwithin{equation}{section}

\begin{document}

\title[]
{A remark on approximation on totally real sets.}

\author[]{ Bo Berndtsson}

\address{B Berndtsson :Department of Mathematics\\Chalmers University
  of Technology 
  and the University of G\"oteborg\\S-412 96 G\"OTEBORG\\SWEDEN,\\}

\email{ bob@math.chalmers.se}

\begin{abstract}
{We give a new proof of  a theorem  on approximation of continuous
  functions on totally real sets.}
\end{abstract}

\bigskip

\maketitle

\section{Introduction}

Let $\Omega$ be a pseudoconvex  open set in $\C^n$, and let $\phi$ be
a  $C^2$-smooth nonnegative plurisubharmonic function in $\Omega$, satisfying
$i\ddbar\phi\geq\delta\beta$, where $\beta$ is the Euclidean volume
form and $\delta>0$. Then
$$
E:=\{\phi =0\}
$$
is a {\it totally real} set. 

Associated to the function $\phi$ and to any positive number $k$ we
have the orthogonal projection operator $P_k$ from 
$$
L^2(\Omega,e^{-k\phi})
$$
to
$$
A^2(\Omega,e^{-k\phi}),
$$
the latter space being the Bergman space, i e the subspace of
holomorphic functions in $L^2(\Omega,e^{-k\phi})$. 

We shall prove the following theorem:
\begin{thm}
Let $u$ be a smooth function of compact support in $\Omega$. Let
$K$ be a compact subset of $\Omega$. Then
$$
\sup_{E\cap K} |u-P_k(u)|\leq C/\sqrt k,
$$
for some constant $C$.

\end{thm}

In particular it follows that any continuous function on $E$ can be
approximated uniformly on compacts of $E$ by functions holomorphic in
$\Omega$. Since any totally real submanifold of class $C^1$ can be
given as the zeroset of a strictly plurisubharmonic function defined in some
neighbourhood of the manifold, Theorem 1.1 contains the theorem of
Hörmander and Wermer, \cite{H-W}and Nirenberg-Wells, \cite{N-W}, as
well as the generalization to the $C^1$-case of Harvey and Wells,
\cite{Harvey-Wells}.  

Just like in the original proofs in \cite{H-W} and \cite{N-W}, the proof of
Theorem 1.1 is based on Hörmander's $L^2$-estimates for the
$\dbar$-operator. One difference between the proofs is that we will
use the weight factor $e^{-k\phi}$ in the estimates. We shall then
use the following consequence of Hörmander's theorem.

\begin{thm} Let $\phi$ be a  plurisubharmonic function in a
  pseudoconvex domain $\Omega$, satisfying $i\ddbar\phi\geq
  \delta\beta$.
Let $f$ be a $\dbar$-closed form of bidegree $(0,1)$ in $\Omega$, and
let $v$ be the $L^2$-minimal solution to the equation $\dbar v=f$ in 
$L^2(\Omega,e^{-k\phi})$. Then, for $k>0$,
$$
\int |v|^2 e^{-k\phi}\leq C/k \int |f|^2 e^{-k\phi}.
$$
\end{thm}
The important feature of this theorem here is that the estimates gets
better as $k$ increases. The same estimate holds if we replace $k\phi$
by $\phi_k$ where $i\ddbar\phi_k\geq k\delta\beta$. 

The next step
in the proof is the observation that something similar happens in
uniform norms, at least if we shrink the domain a little. This is not
entirely trivial, but nor is it
a deep observation - the shrinking of the domain avoids the main
difficulty in passing from $L^2$ to uniform norms. The main point in
the proof is a variant of the Donnelly-Fefferman trick.

\begin{thm} Let $\phi$ be a  plurisubharmonic function in a
  pseudoconvex domain $\Omega$, satisfying $i\ddbar\phi\geq
  \delta\beta$. Let $v$ be the $L^2(\Omega,e^{-k\phi})$-minimal solution
  to the equation $\dbar v=f$. Then, if $K$ is a compact subset of $\Omega$ 
$$
\sup_{K}|v|^2  e^{-k\phi}\leq \frac{C_{\delta,K}}{k} \sup_{\Omega} |f|^2
e^{-k\phi}.
$$
\end{thm}

We apply  Theorem 1.3 to $f=\dbar u$ where $u$ is, say, a test function.
Since
$$
u-P_k(u)= v,
$$
it follows that with $E=\phi^{-1}(0)$
$$
\sup_{E\cap K}|u-P_k(u)|^2\leq \frac{C_{\delta,K}}{k} \sup_{\Omega} |f|^2
e^{-k\phi}\leq \frac{C_{\delta , K}}{k} \sup_{\Omega} |f|^2,
$$ 
if $\phi\geq 0$, and Theorem 1.1 follows. 

It is interesting to note that the proof of the  approximation theorem
here also has some
features in common with the proof of (a generalization of) the
Hörmander-Wermer theorem of Baouendi and Treves, \cite{B-T}. Their
proof is based on convolution with a Gaussian kernel, whereas here we
apply the Bergman kernel. However, in the model case $\phi= x^2$, the
Bergman kernel is Gaussian, so the two proofs are actually quite 
similar in this case.

If $K$ is a compact subset of $E$, which is moreover polynomially
convex,  we can choose the weight function
$\phi$ so that it has logarithmic growth at infinity. The holomorphic
functions $P_k(u)$ are then polynomials of degree $k$, and Theorem 1.1
 estimates the degree of approximation by $1/\sqrt k$, if $u$ is
of class $C^1$. This is not quite as good as one would expect; at
least if $E$ is a smooth manifold the right degree of approximation
should be 
$1/k$. Possibly this flaw comes from letting $E$  be a quite general 
set. At any rate it seems hard to do better with the methods in this
paper. 

One might also notice that  Theorem 1.3 works equally well in
unbounded domains, so minor modifications should give Carleman-Type
approximation as well, see \cite{Carleman}, \cite{Manne}.

Finally, I would like to thank  Said Asserda and the referee for pointing out
several inaccuracies and obscurities in the first version of this paper.

\section{Proof of Theorem 1.3}

Theorem 1.3 is not really new - it follows readily from  results in
\cite{B1} and \cite{B2} and in particular \cite{Delin} - but here we
shall indicate a concise proof 
for the case at hand. The main ingredient is a variant of the
Donnelly-Fefferman trick which will give us an Agmon-type estimate.
\begin{thm} Assume $i\ddbar\phi\geq 5\beta$. Then, with notation as in
  the introduction, 
  for any $a$ in $\C^n$
\begin{equation}
\int |v|^2 e^{-k\phi -\sqrt k|z-a|}\leq C/k \int |f|^2 e^{-k\phi
  -\sqrt k|z-a|}. 
\end{equation}
\end{thm}

\begin{proof}
Assume for simplicity of notation that $a=0$. Let $\psi$ be the convex
function defined by
$$
\psi(t)=t
$$
for $t\geq 1$ and 
$$
\psi(t)=t^2/2 +1/2
$$
for $t<1$. Let $\chi(z)=\psi(|z|)$. Then $|\chi -|z||$ is bounded,
so it is enough to prove 2.1 with $\sqrt k |z|$ replaced by
$\chi(z\sqrt k)$. Moreover $\partial\chi$ is bounded by 1, and
$i\ddbar\chi\leq\beta$. It is 
especially the last property that is of importance here and is the
reason for introducing the function $\chi$. Put
$$
\chi_k(z)=\chi(z\sqrt k)
$$
and
$$
v_k:= v e^{-\chi_k}.
$$

Since $v$ is orthogonal to all holomorphic functions in
$L^2(\Omega,e^{-k\phi})$, it follows that $v_k$ is orthogonal to all
holomorphic functions for the scalar product in
$L^2(\Omega,e^{-k\phi+\chi_k})$ (we may assume in the
proof that $\Omega$ is bounded so that the $L^2$-spaces do not change
when we vary the weight). Hence $v_k$ is the $L^2$minimal solution to
a certain $\dbar$-equation. Now,
$$
i\ddbar (k\phi-\chi_k)\geq Ck\beta 
$$
so it follows from the Hörmander estimate 1.2  
that
\begin{equation}
\int |v_k|^2 e^{-k\phi+\chi_k}\leq C/k\int |\dbar v_k|^2
e^{-k\phi+\chi_k} .
\end{equation}
The left hand side of equation 2.2 equals

$$ 
\int |v|^2 e^{-k\phi-\chi_k}.
$$
In the right hand side we have
$$
\dbar v_k =(f -v_k\dbar\chi_k)e^{-\chi_k}.
$$
Since $\dbar\chi_k$ is bounded by $\sqrt k$ we can absorb the contribution to
2.2 coming from the second term $v_k\dbar\chi_k$ in the left hand side
of 2.2, and 2.1 follows. 

\end{proof}

We are now ready for the proof of Theorem 1.3. By scaling, we may of course
assume that $ i\ddbar\phi\geq 5\beta$.
Let $a$ be a point in $\Omega_k$ that again for simplicity we  take
equal to 0. After changing frame locally at 0 we can  assume that,
near the origin,
$$
\phi(z)= q(z,\bar z) + o(|z|^2),
$$
 where $q$ is an hermitian form. 

(This means the following: Near the origin we can, since $\phi$ is of
class $C^2$, write
$$
\phi(z)=2\Re P(z) +q(z,\bar z)+o(|z|^2),
$$
where $P$ is a holomorphic polynomial of degree 2. We then write 
$$
v'=v e^{-P(z)}
$$
for $v$ and 
$$
f'=fe^{-P(z)}
$$
for $f$, which changes the weight function $\phi$ to $\phi -2\Re P=q
+o(|z|^2)$.)

\bigskip

\noindent
 Note that $kq(z,\bar z)$ is bounded by a constant when
 $|z|<1/k$. This constant certainly depends on the point $a$ that we
 have taken equal to 0, but it is
 uniform as long as $a$ ranges over a compact subset of $\Omega$. We then get
$$
\int_{|z|^2<1/k}|v|^2\leq C/k\sup |f|^2 e^{-k\phi}\int e^{-\sqrt
  k|z|}\leq C\sup |f|^2 e^{-k\phi}/k^{n+1}.
$$

Normalize so that
$$
\sup |f|^2 e^{-k\phi}\leq 1.
$$
Then, in particular, $\dbar v$ is bounded by a constant for
$|z|^2<1/k$. The   inequality  
\be
|v(0)|^2\leq C(k^n\int_{|z|^2<1/k}|v|^2 + \frac{1}{k}\sup_{|z|^2 <1/k} |f|^2 ),
\ee
then shows that
$$
|v(0)|^2\leq 1/k,
$$
which is what we wanted to prove.  To verify 2.3 one can apply the
Bochner-Martinelli integral formula
$$
v(0)=c_n\int (\dbar (v \xi(|z|^2/k)\cdot \partial |z|^{2n-2}),
$$
where $\xi(t)$ is a smooth function that equals 1 for $t<1/2$ and 0
for $t>1$.

\bigskip

\def\listing#1#2#3{{\sc #1}:\ {\it #2}, \ #3.}

\end{document}